\documentclass[12pt,a4paper]{amsart}
\usepackage{mathtools}
\usepackage{amssymb}
\usepackage{amscd}
\usepackage[mathscr]{eucal}
\usepackage{latexsym}

\usepackage{verbatim}
\usepackage{setspace}
\usepackage{geometry}
\usepackage{float}
\usepackage[english]{babel}
\usepackage[utf8]{inputenc}

\usepackage{graphicx}
\usepackage{subcaption}
\usepackage[dvipsnames]{xcolor}
\usepackage{framed}
\usepackage[all]{xy}
\usepackage{tikz}
\usetikzlibrary{trees}

\usepackage{enumitem}
\usepackage{listings}

\usepackage{hyperref}
\usepackage{cleveref}

\theoremstyle{plain}
\newtheorem{thm}{Theorem}[section]
\newtheorem{lem}[thm]{Lemma}

\theoremstyle{definition}
\newtheorem{example}[thm]{Example}
\newtheorem{defi}[thm]{Definition}
\newtheorem{rem}[thm]{Remark}

\def\qed{\hfill $\Box$}

\definecolor{supergray}{gray}{0.95}
\definecolor{hgray}{gray}{0.9}
\DeclareMathOperator{\MOD}{\rm{MOD}}
\DeclareMathOperator{\Int}{\rm{Int}}
\DeclareMathOperator{\St}{\rm{St}}
\Crefname{figure}{Figure}{Figures}

\title[Non-singular extensions of circle-valued Morse functions]{Non-singular extensions of circle-valued Morse functions}
\author[KOKI IWAKURA]{Koki Iwakura}
\address{Joint Graduate School of Mathematics for Innovation, Kyushu University, Motooka 744, Nishiku, Fukuoka 819-0395, Japan.}
\email{iwakura.koki0105@gmail.com}

\begin{document}

\maketitle
\begin{abstract}
In this paper, we consider the non-singular extension problem for circle-valued Morse functions on closed orientable surfaces.
The problem asks, given a circle-valued Morse function $f\colon M\to S^{1}$ on a closed orientable surface $M$, under what condition there exist a compact orientable 3-dimensional manifold $N$ with $\partial N = M$ and a submersion $G\colon N \to S^{1}$ such that $G|_{\partial N}=f$. 
We provide necessary and sufficient conditions for the existence of a non-singular extension of a circle-valued Morse function as the main theorem when a submersion on a collar neighborhood is given.
\end{abstract}

\section{Introduction}
In this paper, we consider the ``non-singular extension problem" for circle-valued Morse functions.
Let $M$ be a closed orientable (possibly disconnected) surface and $g$ be a submersion of $M\times[0,1)$ into $S^{1}$ such that $g|_{M\times\{0\}}$ is a circle-valued Morse function. 
Then, under what condition does there exist a compact orientable $3$-dimensional manifold $N$ with boundary $\partial N=M$ and a submersion $G\colon N\to S^{1}$ that coincides with $g$ on the collar neighborhood?
The answer to this problem is given in the following theorem.

\begin{thm}
Let $M$ be a closed orientable surface and $g\colon M\times[0,1)\to S^{1}$ be a submersion such that $f\coloneq g|_{M\times\{0\}}$ is a circle-valued Morse function. 
Then, there exist a compact orientable $3$-dimensional manifold $N$ and a non-singular extension $G\colon N\to S^{1}$ of $g$ if and only if there exist a finite disjoint union $V$ of finite graphs and some circles, a continuous map $h\colon V\to S^{1}$ which is an immersion on the edges and circle components of $V$, and an allowable collapse $C\colon W_{f}^{\pm}\to V$ for $\overline{f}$ and $h$.
\end{thm}

The notion of a ``collapse" was introduced by Curley in \cite{C} for real-valued Morse functions on closed orientable surfaces. 
A ``collapse" contains information about the relationships between the map on the surface and that on the $3$-dimensional manifold.
Furthermore, Curley introduced the notion of an ``allowable collapse" to obtain a certain correspondence between two critical points of a real-valued Morse function that increase and decrease the genera of level surfaces. 
In this paper, we modify the definitions of a ``collapse" and an ``allowable collapse" to make them compatible with circle-valued Morse functions and prove the above theorem.
For these functions, the difficulty arises because, unlike two points in $\mathbb{R}$, we cannot impose an order on two points in $S^{1}$.
It should be noted that this work is a direct generalization of the main theorem in \cite{C}.
\par
The ``non-singular extension problem" is classical, and here, let us mention some previous studies.
The non-singular extension problem for real-valued Morse functions on $1$-dimensional manifold was studied by Blank--Laudenbach in \cite{BL}. 
The non-singular extensions of Morse functions on closed orientable surfaces were studied by Curley in \cite{C}. 
Moreover, Curley considered the non-singular extension problem of a Morse function on $S^{2}$ to a submersion on $D^{3}$. 
In \cite{A}, Arnold formulated the non-singular extension problem.  
Barannikov \cite{B} studied the non-singular extension problem for Morse functions on spheres of arbitrary dimension by using Morse complexes and obtained necessary conditions in terms of combinatorial properties of a Morse complex.
Furthermore, Seigneur \cite{S} studied the same problem and obtained the necessary conditions in terms of algebraic properties of a Morse complex.
Laroche \cite{L} studied the non-singular extension problem for Morse functions on closed non-orientable surfaces and obtained necessary and sufficient conditions.
Furthermore, Laroche considered the non-singular extension problem of a Morse function on a Klein bottle to a solid Klein bottle. 
Given a Morse function on a closed orientable (possibly disconnected) surface, Iwamoto \cite{I} succeeded in extending it to a submersion on a collar appropriately so that it can be further extended to a submersion of a compact orientable connected $3$-dimensional manifold by using the Curley's result in \cite{C}.
\par
The present paper is divided into four parts.
In Section 2, we introduce some notations and lemmas used throughout the paper. 
In Section 3, we give some concrete examples, containing those that arise only when we consider circle-valued Morse functions. 
In Section 4, we prove the main theorem.
Throughout the paper, manifolds and maps between them are assumed to be of class $C^{\infty}$ unless otherwise specified.

\section{Preliminaries}
In this section, we prepare some necessary materials used in this paper. 
Note that we assume that the circle as the target manifold of circle-valued Morse functions is oriented for the rest of this paper. 

\begin{defi}[Circle-valued Morse function]
Let $M$ be a closed surface. 
Then a smooth map $f\colon M\to S^{1}$ is called a \emph{circle-valued Morse function} if its critical points are all non-degenerate and the critical values are all distinct.
\end{defi}

\begin{rem}\label{rem:1}
In the present paper, we allow the empty set as a critical point set of circle-valued Morse functions.
This means that a circle-valued Morse function can be a submersion.
\end{rem}

The following lemma determines local forms of circle-valued Morse functions around critical points.
This is obtained by Morse's lemma in \cite{Mi}.
\begin{lem}\label{lem:1}
Let $M$ be a closed surface and $f\colon M\to S^{1}$ be a circle-valued Morse function. 
Let $p\in M$ be a critical point of $f$.
Then, there exist local coordinates $(x,y)$ around $p\in M$ and $z$ around $f(p)\in S^{1}$ such that $f$ is of the form
\begin{equation*}
(x,y)\mapsto z=\pm x^{2}\pm y^{2},
\end{equation*}
where the local coordinate $z$ of $S^{1}$ is consistent with the orientation.
\end{lem}

Let us give local normal forms of submersions around critical points of the circle-valued Morse function on the boundary.
We omit the proof of the following lemma since it is similar to that of \cite[Proposition 1]{Shi}.

\begin{lem}\label{lem:2}
Let $N$ be a compact orientable $3$-dimensional manifold and $G\colon N\to S^{1}$ be a submersion such that $G|_{\partial N}$ is a circle-valued Morse function. 
Let $p \in \partial N$ be a critical point of $G|_{\partial N}$.
Then, there exist local coordinates $(X,Y,Z)$ around $p\in N$ and $W$ around $G(p)\in S^{1}$ such that $G$ is of the form
\begin{equation*}
(X, Y, Z)\mapsto W=\pm X^{2}\pm Y^{2}\pm Z,
\end{equation*}
where $\{Z=0\}$ corresponds to the boundary of $N$, $\{Z>0\}$ corresponds to the interior of $N$, and the local coordinate $W$ of $S^{1}$ is consistent with the orientation.
\end{lem}

Let us define an equivalence relation on compact manifolds induced by a smooth map from those manifolds into $S^{1}$. 
We define the ``Reeb graph" of a map following \cite[Definition 2.1]{G} as follows.

\begin{defi}[Reeb graph]\label{def:1}
Let $M$ be a compact manifold and $f\colon M\to S^{1}$ be a smooth map. 
A relation ``$\sim$" on $M$ is defined as follows:
for $x,y\in M$, $x\sim y$ if $f(x)=f(y)$ and  $x,y$ are in the same connected component of $f^{-1}(f(x))$.
This is an equivalence relation on $M$ and the quotient space $W_{f}\coloneqq M/\sim$ is called the \emph{Reeb graph} of $f$. Furthermore, it is easy to prove that there is a unique continuous map $\overline{f}$ that makes the following diagram commutative:
\[
  \xymatrix{
    M \ar[r]^{f} \ar[d]_{q_{f}} & S^{1} \\
    W_{f} \ar[ru]_{\overline{f}},
  }
\]
where the map $q_{f}$ is the quotient map.
\end{defi}

We define the vertices of the Reeb graph of a circle-valued Morse function on a closed manifold as corresponding to the level sets that contain its critical points.
According to \cite[Proposition 3.1]{RL}, the Reeb graph of such a function is a finite disjoint union of finite graphs and circles. 
Since we also allow submersions as circle-valued Morse functions, as noted in \cref{rem:1}, Reeb graphs may have circle components.
Furthermore, let us define the vertices of the Reeb graph of a submersion on a manifold with non-empty boundaries, where the submersion is a circle-valued Morse function on the boundaries.
In this case, we take vertices corresponding to the level sets that contain the critical points of the function on the boundaries.
The Reeb graph is also a finite disjoint union of finite graphs and circles, as follows from \textrm{Lemma~2.4}.

\begin{defi}[Non-singular extension]\label{def:3}
Let $M$ be a closed orientable surface and $g\colon M\times[0,1)\to S^{1}$ be a submersion such that the restriction to $M\times\{0\}$ is a circle-valued Morse function. 
We assume that there exist compact orientable $3$-dimensional manifold $N$ and a submersion $G\colon N\to S^{1}$ such that each connected component of $N$ has non-empty boundary, $\partial N=M$, and the diagram
\[
  \xymatrix{
    M\times[0,1) \ar[r]^{~~~~g} \ar[d]_{i\;} & S^{1} \\
    N \ar[ru]_{G}
  }
\]
is commutative for some $i$, which is a collar identifying $M\times\{0\}$ with $\partial N$. 
Then, we call $G$ a \emph{non-singular extension} of $g$.
\end{defi}

Let $M$ be a closed orientable surface and $g\colon M\times [0,1)\to S^{1}$ be a submersion such that $f\coloneq g|_{M\times\{0\}}$ is a circle-valued Morse function. 
Let $p\in M\times\{0\}$ be a critical point of $f$ and $w\in T_{p}(M\times[0,1))$ be an outward vector. 
Note that $dg_{p}(w)$ is not zero since $g$ is a submersion.
If $dg_{p}(w)$ fits the orientation of $S^{1}$, then we label the vertex $q_{f}(p)$ of the Reeb graph $W_{f}$ of $f$ by ``$+$", otherwise ``$-$".
Note that the labels do not depend on the choice of the outward vectors. 

\begin{defi}[Labeled Reeb graph]
The Reeb graph of $f$, in which each vertex is labeled by ``$+$" or ``$-$" according to the above rules, is called a \emph{labeled Reeb graph} of $f$ and is denoted by $W_{f}^{\pm}$. 
\end{defi}

\begin{defi}[Collapse]\label{def:5}
Let $M$ be a closed orientable surface and $g\colon M\times[0,1)\to S^{1}$ be a submersion such that $f\coloneq g|_{M\times\{0\}}$ is a circle-valued Morse function. 
Let $V$ be a finite disjoint union of finite graphs and some circles, and let $h\colon V\to S^{1}$ be a continuous map that is an immersion on each edge and circle, where edges are identified with closed intervals.
Then, a continuous map $C\colon W_{f}^{\pm}\to V$ is called a \emph{collapse} if the following conditions are satisfied:

\begin{figure}[t]
\centering
        \includegraphics[width=98mm]{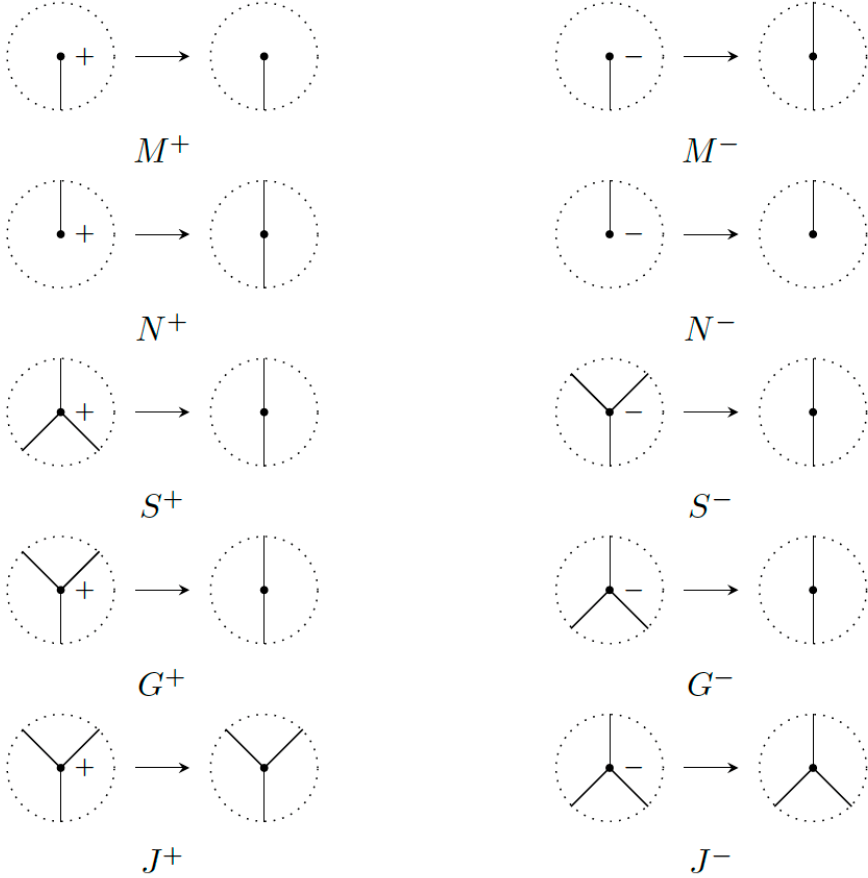}
        \vspace{10pt}
\caption{
List of local behaviors of a collapse around the vertices of $W_{f}^{\pm}$.
The signs on the symbols $M$, $N$, $S$, $G$, and $J$ correspond to the signs assigned to the vertices of $W_{f}^{\pm}$ and those symbols are adopted from \cite{C}.
}
\vspace{0pt}
\label{fig:2}
\end{figure}

\begin{enumerate}
\setlength{\parskip}{0cm}
  \setlength{\itemsep}{0cm}
\item The map $C$ maps the vertices of $W_{f}^{\pm}$ to those of $V$, bijectively.
\item The map $C$ maps the circle components of $W_{f}^{\pm}$ to those of $V$, surjectively.
\item The diagram
\[
  \xymatrix{
    W_{f}^{\pm} \ar[r]^{\overline{f}} \ar[d]_{C} & S^{1} \\
    V \ar[ru]_{h}
  }
\]
is commutative: i.e., we have $\bar{f}=h\circ C$.
\item Around each vertex of $W_{f}^{\pm}$, the map $C$ is locally equivalent to one of the maps as depicted in \Cref{fig:2}.
\end{enumerate}
\end{defi}

Let us define a property for a collapse. 
To state the property, let us introduce the relation ``$\preceq$" for the vertices of $V$ as follows. 
For vertices $x,y\in V$, we write $x\preceq y$ if there exists a path $\sigma\colon [0,1]\to V$ such that 
\begin{enumerate}
\item $\sigma(0)=x$ and $\sigma(1)=y$, 
\item the map $h\circ\sigma\colon [0,1]\to S^{1}$ is an orientation preserving immersion. 
\end{enumerate}

\begin{defi}[Allowable]
In the context of \Cref{def:5}, let $\mathscr{G}^{+}$ (resp. $\mathscr{G}^{-}$) be the subset of the vertex set of  $W_{f}^{\pm}$ consisting of the vertices corresponding to $G^{+}$ (resp. $G^{-}$) as described in \Cref{fig:2}. 
A collapse $C$ is called \emph{allowable} if there exists a bijection $\gamma\colon\mathscr{G}^{+}\to \mathscr{G}^{-}$ such that $C(\gamma(a))\preceq C(a)$ for every $a\in\mathscr{G}^{+}$.
\end{defi}

\section{Examples}

\begin{example}
In \textrm{Figure~2(a)}, the figure on the left depicts a submersion $g_{1}\colon S^{2}\times[0,1)\to S^{1}$, while that on the right depicts the labeled Reeb graph of the circle-valued Morse function $f_{1}$ on its boundary. 
In \textrm{Figure~2(b)}, the submersion $G_{1}\colon D^{3}\to S^{1}$ on the left is a non-singular extension of $g_{1}$.
The figure on the right shows the Reeb graph of $G_{1}$.
\begin{figure}[b]
    \centering
    \begin{subfigure}[b]{0.7\textwidth} 
        \centering
        \includegraphics[width=85mm]{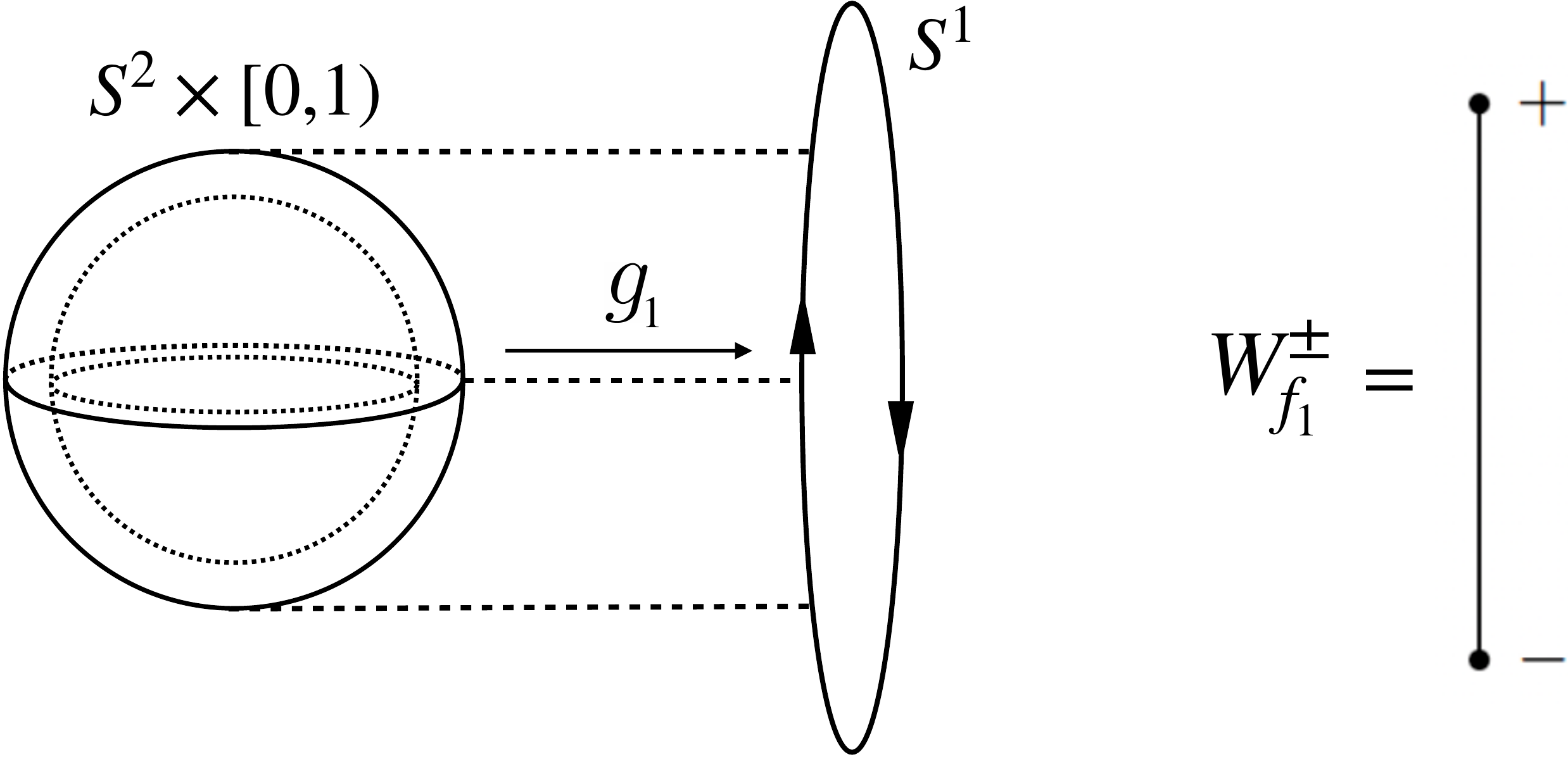}
        \vspace{-3pt}
        \caption{$g_{1}$ and $W^{\pm}_{f_{1}}$.}
        \label{fig:3a}
    \end{subfigure}\\[2ex]
            \vspace{-7pt}
    \begin{subfigure}[b]{0.7\textwidth} 
        \centering
        \includegraphics[width=85mm]{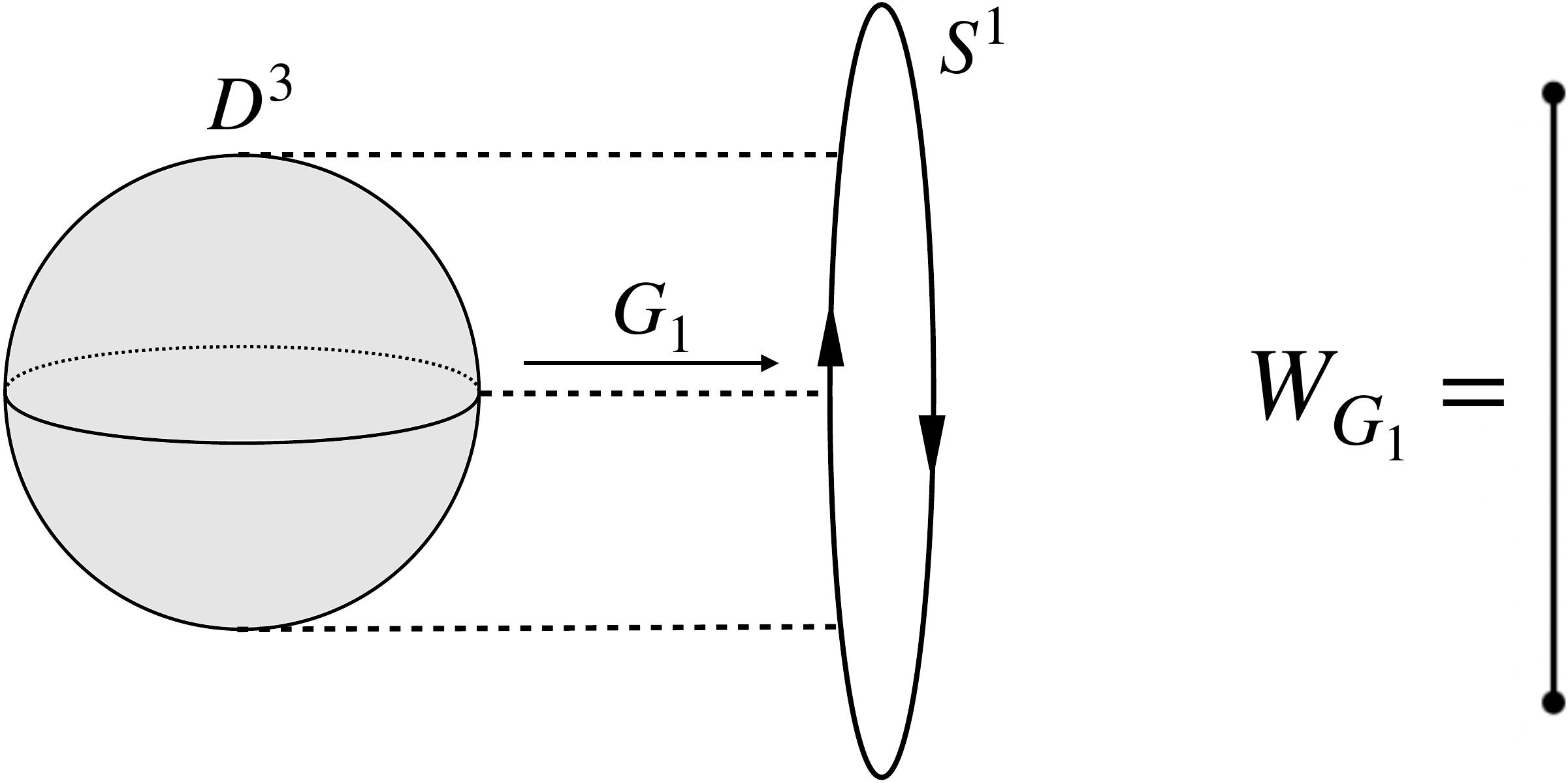}
                \vspace{-3pt}
        \caption{$G_{1}$ and $W_{G_{1}}$.}
        \label{fig:3b}
    \end{subfigure}
        \vspace{-11pt}
    \caption{}
    \label{fig:3}
\end{figure}
\end{example}

\begin{example}
A submersion $g_{2}\colon (S^{2}\sqcup S^{2})\times[0,1)\to S^{1}$ is depicted on the left of \textrm{Figure~3(a)}, and the corresponding labeled Reeb graph of $f_{2}\coloneq g_{2}|_{(S^{2}\sqcup S^{2})\times\{0\}}$ is shown on the right. 
In \textrm{Figure~3(b)}, the submersion $G_{2}\colon N_{2}\to S^{1}$ on the left is a non-singular extension of $g_{2}$, where the $3$-dimensional manifold $N_{2}$ is obtained by removing a $3$-dimensional open ball from $D^{3}$. 
The Reeb graph of $G_{2}$ is shown on the right.
\begin{figure}[t]
    \centering
    \begin{subfigure}[b]{0.7\textwidth} 
        \centering
        \includegraphics[width=85mm]{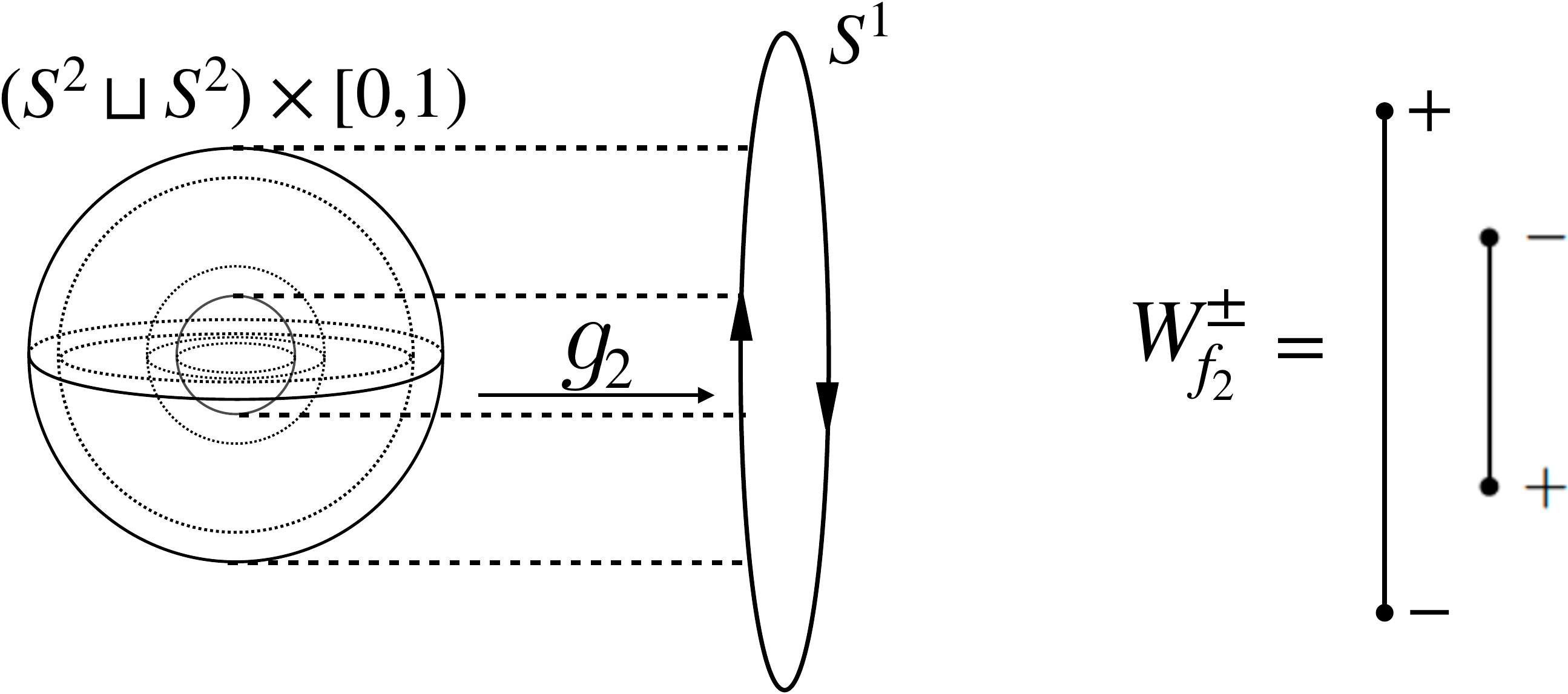}
         \vspace{-3pt}
        \caption{$g_{2}$ and $W^{\pm}_{f_{2}}$.}
        \label{fig:4a}
    \end{subfigure}\\[2ex]
    \vspace{-7pt}
    \begin{subfigure}[b]{0.7\textwidth} 
        \centering
        \includegraphics[width=85mm]{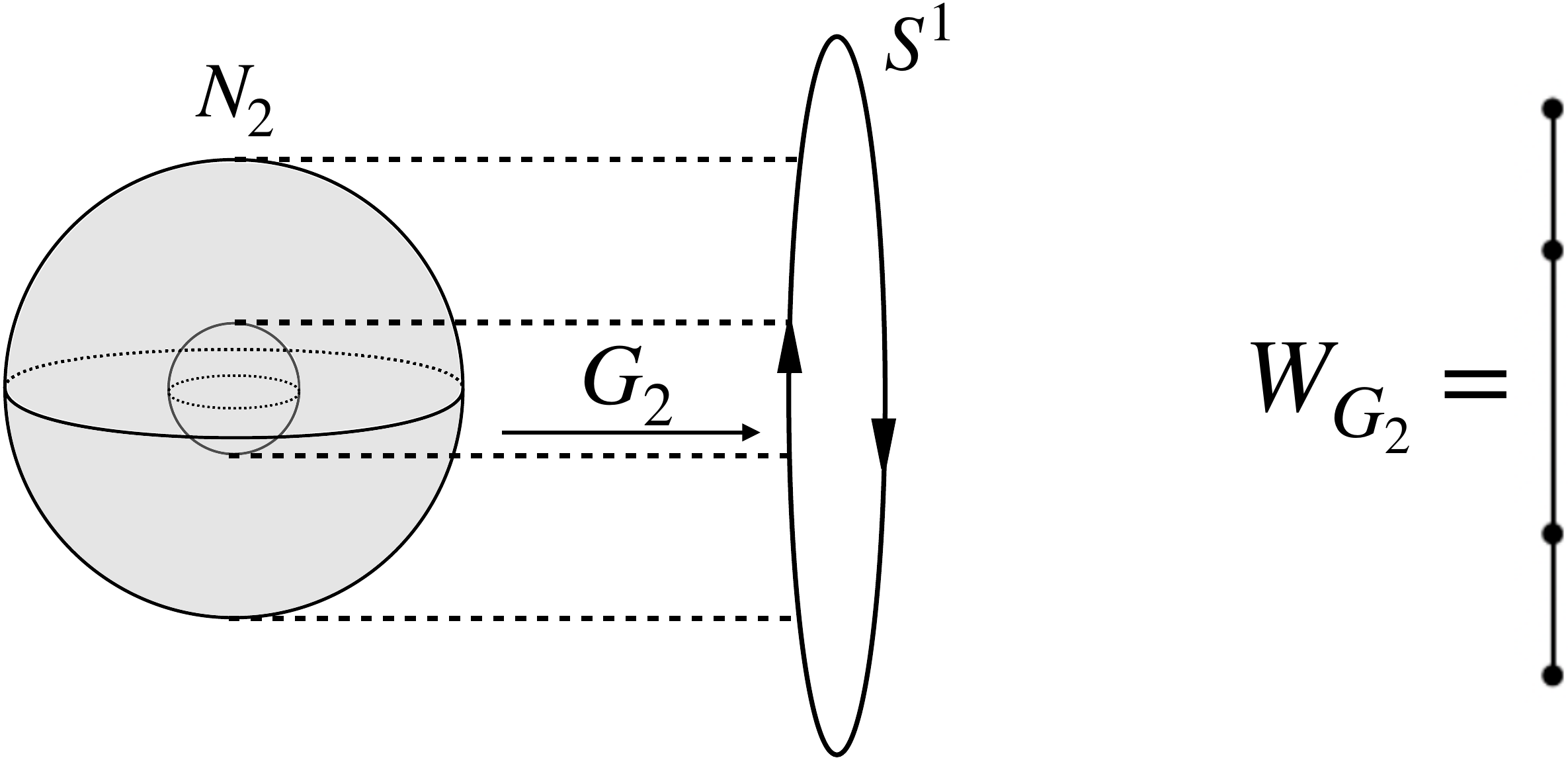}
         \vspace{-3pt}
        \caption{$G_{2}$ and $W_{G_{2}}$.}
        \label{fig:4b}
    \end{subfigure}
    \vspace{-11pt}
    \caption{}
    \label{fig:4}
\end{figure}
\begin{figure}[t]
    \centering
    \begin{subfigure}[b]{0.7\textwidth} 
        \centering
        \includegraphics[width=85mm]{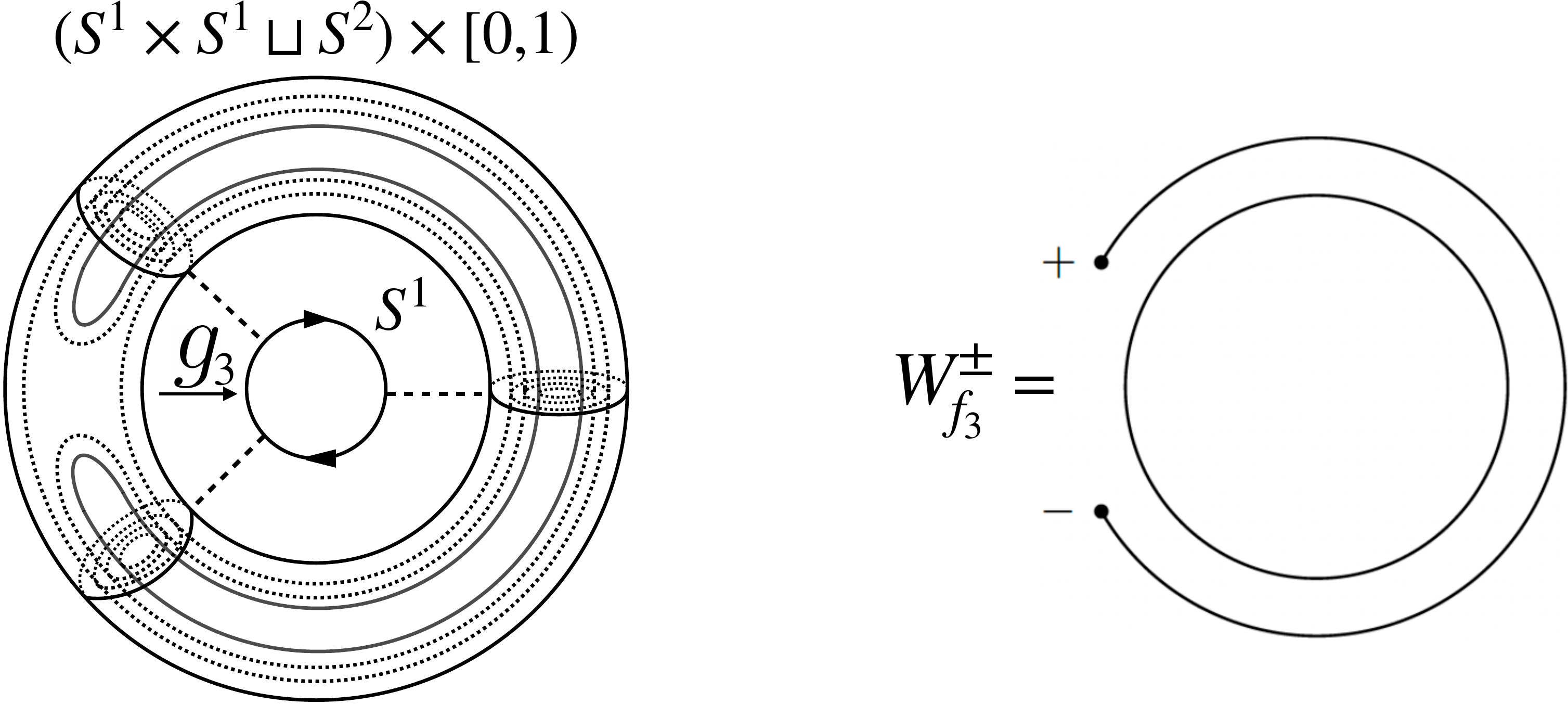}
                \vspace{-5pt}
        \caption{$g_{3}$ and $W^{\pm}_{f_{3}}$.}
        \label{fig:6a}
    \end{subfigure}\\[2ex]
            \vspace{-5pt}
    \begin{subfigure}[b]{0.7\textwidth} 
        \centering
        \includegraphics[width=85mm]{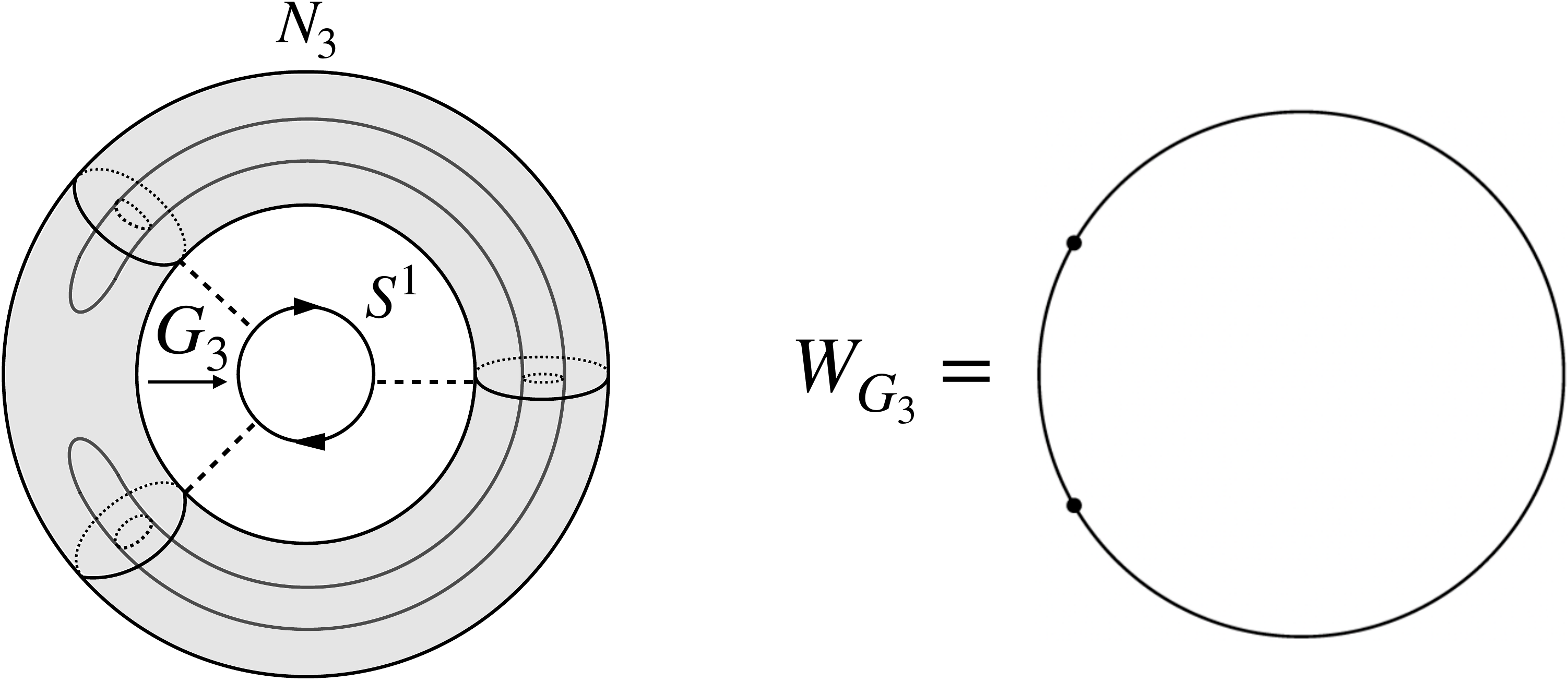}
                \vspace{-5pt}
        \caption{$G_{3}$ and $W_{G_{3}}$.}
        \label{fig:6b}
    \end{subfigure}
        \vspace{-8pt}
    \caption{}
    \label{fig:6}
\end{figure}
\end{example}

\begin{example}
A submersion $g_{3}\colon (S^{1}\times S^{1}\sqcup S^{2})\times[0,1)\to S^{1}$ and the labeled Reeb graph of $f_{3}\coloneq g_{3}|_{(S^{1}\times S^{1}\sqcup S^{2})\times\{0\}}$ are depicted in \textrm{Figure~4(a)}. 
In \textrm{Figure~4(b)}, the left figure depicts the submersion $G_{3}\colon N_{3}\to S^{1}$, which is a non-singular extension of $g_{3}$.
The $3$-dimensional manifold $N_{3}$ is obtained by removing a $3$-dimensional open ball from $S^{1}\times D^{2}$.
The Reeb graph of $G_{3}$ is depicted on the right. 
\end{example}

\section{Proof of \textrm{Theorem~1.1}}
\par
Let us assume that a compact orientable $3$-dimensional manifold $N$ and a non-singular extension $G\colon N\to S^{1}$ of $g$ are given. 
Then, we obtain the commutative diagram.
\[
  \xymatrix{
    W_{f}^{\pm} \ar[r]^{\overline{f}} \ar[d]_{C} & S^{1} \\
    W_{G} \ar[ru]_{\overline{G}}
  }
\]
by the commutative diagram in \Cref{def:3}, where the map $C$ is induced by the inclusion $i|_{M\times\{0\}}$.
By \textrm{Lemma~2.4}, the Reeb graph $W_{G}$ of $G$ is a finite disjoint union of finite graphs and some circles, where the vertices correspond bijctively to those of $W_{f}^{\pm}$ and the circle components are derived from those of $W_{f}^{\pm}$.
Since the local forms of $G$ around the critical points of $f$ are determined by \Cref{lem:2}, the behaviors of $C$ around the vertices of $W_{f}^{\pm}$ are as listed in \textrm{Figure~1}. 
Therefore, the map $C$ is a collapse.

\begin{figure}[b]
\centering
\includegraphics[width=140mm]{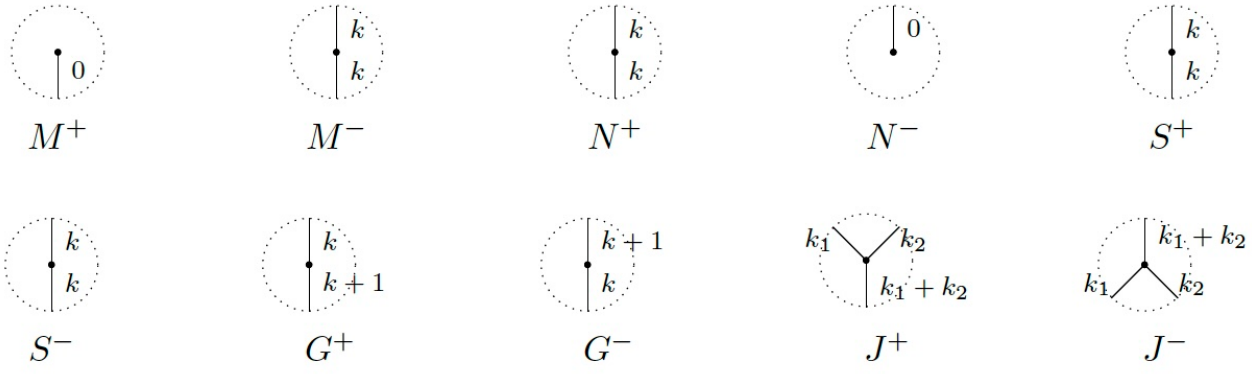}
\caption{
We depict the neighborhoods of the vertices of $W_{G}$, where the top and bottom portions of the dotted circles are aligned with $\overline{G}$ according to the orientation of $S^{1}$.
The numbers on the edges indicate the genera of the level surfaces of $G$ corresponding to those edges.
The symbols correspond to those in \Cref{fig:2}.
}
\label{fig:7}
\vspace{-4pt}
\end{figure}

Let us prove that the collapse $C$ is allowable. 
In \textrm{Figure~5}, since the topological types of level surfaces of $G$ change only at the critical points of $f$, the numbers on the edges, which represent their genera, are constant. 
For each case in \textrm{Figure~5}, changes in the numbers at the vertices can be understood by analyzing the Euler characteristic of the level surfaces. 
We take paths associated with these numbers, where the direction of each path is opposite to the orientation of $S^{1}$ when mapped by $h$. 
These paths originate at $G^{+}$ and terminate at $G^{-}$, with merging at $J^{+}$ and branching at $J^{-}$ according to the numbers on the edges.
The generating and terminating of paths occur only at $G^{+}$ and $G^{-}$, and this happens without excess and deficiency, as the total number of paths remains conserved.   
Therefore, there exist paths from each $G^{-}$ vertex to each $G^{+}$ vertex bijectively, and their images by $h$ are orientation preserving immersions.
Hence, we obtain a bijection $\gamma\colon\mathscr{G}^{+}\to\mathscr{G}^{-}$ that holds $C(\gamma(a))\preceq C(a)$ for any $a\in\mathscr{G}^{+}$.
Thus, $C$ is an allowable collapse.

Let us prove the converse. 
Denote the vertices of $W_{f}^{\pm}$ by $v_{j}$, those of $V$ by $w_{j}\coloneq C(v_{j})$, and $\overline{f}(v_{j})$ by $a_{j}$ $(j=1,2,\ldots,n)$, where the subscripts are assigned so that $a_{1},a_{2},\ldots,a_{n}$ are ordered sequentially in reverse to the orientation of $S^{1}$. 
Select $b_{j}\in S^{1}\setminus\{a_{1},a_{2},\dots, a_{n}\}$ $(j=1,2,\ldots,n+1)$ so that they are arranged in the exact order $b_{1}, a_{1}, b_{2}, \ldots, a_{n-1}, b_{n}, a_{n}, b_{n+1}$, sequentially in the direction opposite to the orientation of $S^{1}$ as depicted in \textrm{Figure~6}.

\begin{figure}[t]
\vspace{-3pt}
\centering
\includegraphics[width=55mm]{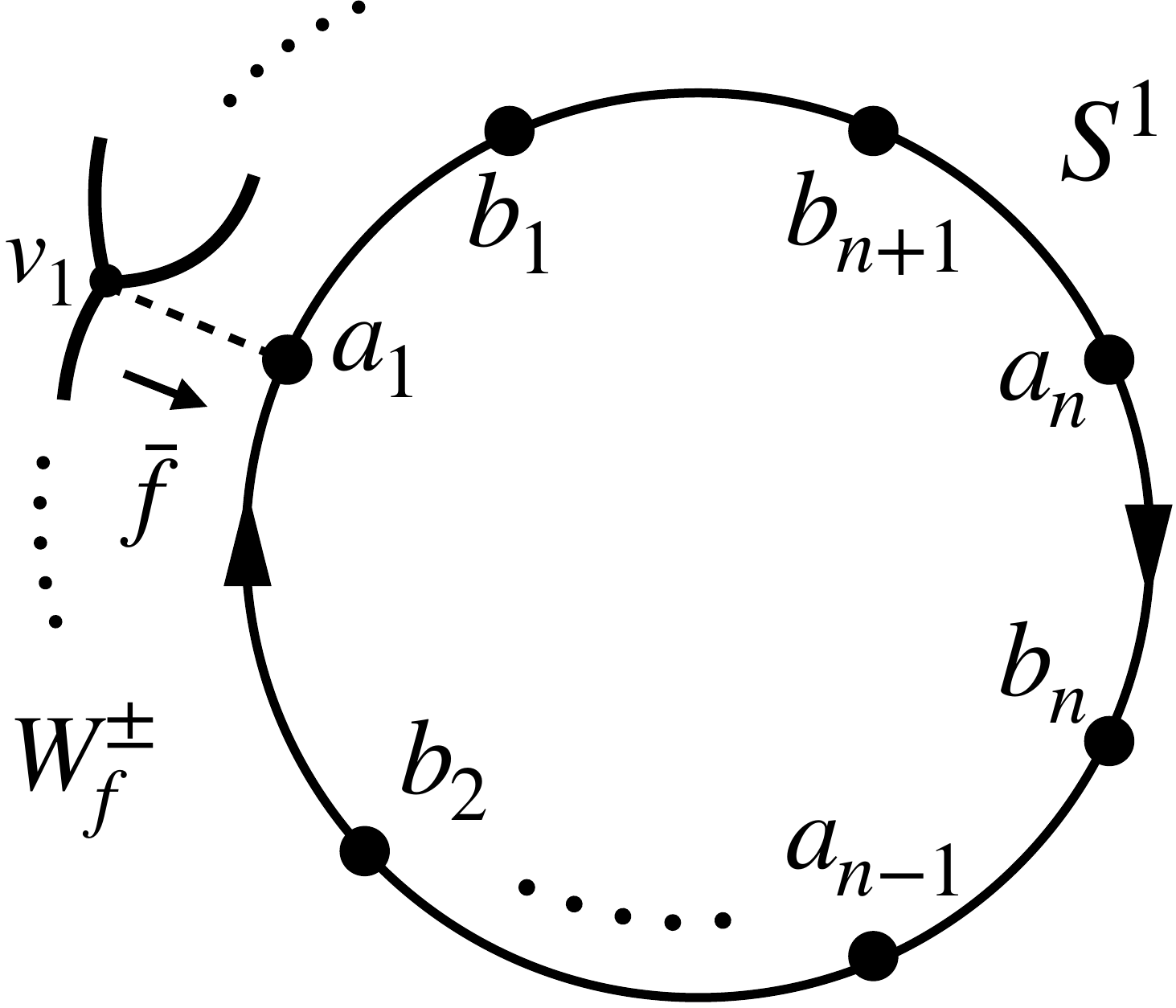}
\vspace{5pt}
\caption{
We depicts the summary figure of the arrangements of $a_{1},a_{2}, \ldots, a_{n}$, $b_{1},b_{2}, \ldots, b_{n+1}$ on $S^1$.
}
\vspace{0pt}
\end{figure}

First of all, we take surfaces corresponding to $b_{1}$ via $h$ as follows.
The number of the connected components equals the number of the edges of $V$ that intersect with $h^{-1}(b_{1})$.
Assuming $C$ is allowable, we assign genera as the number of the paths on the edges, and the numbers follow \textrm{Figure~5} since $h\circ\sigma$ is an orientation preserving immersion for each path $\sigma$. 
The boundary components of each surface correspond to the elements of the inverse image under $C$ of the point in $V$ corresponding to $b_{1}$ via $h$.
We denote these surfaces by $L_{1}$.

We construct a $3$-dimensional manifold $N$ and a non-singular extension $G\colon N\to S^{1}$ of $g$ from a $3$-dimensional manifold $N_{n}$ and a submersion $G_{n}\colon N_{n}\to S^{1}$, which are constructed inductively by using the parts in \textrm{Figure~7} as follows.
At each step of the induction, given a $3$-dimensional manifold $N_{i}$ and a submersion $G_{i}\colon N_{i}\to S^{1}$, we obtain a $3$-dimensional manifold $N_{i+1}$ and a submersion $G_{i+1}\colon N_{i+1}\to S^{1}$ through the processes described later $(i=1,2,\ldots,n-1)$, with the initial condition $N_{0}=L_{1}$.
Here, $N_{i}\subset N_{i+1}$, $G_{i+1}|_{N_{i}}=G_{i}$ $(i=0,1,\ldots,n-1)$, and the boundary components of each connected component of the level surfaces of $G_{i+1}$ correspond to the elements of the inverse image under $C$ of the points in $V$ that correspond via $h$.
In what follows, for a given submersion $G_{i}\colon N_{i}\to S^{1}$, let us denote $G_{i}^{-1}(b_{i+1})\subset N_{i}$ as $L_{i+1}$ $(i=0,1,\ldots,n)$.

\begin{figure}[t]
\centering
\includegraphics[width=80mm]{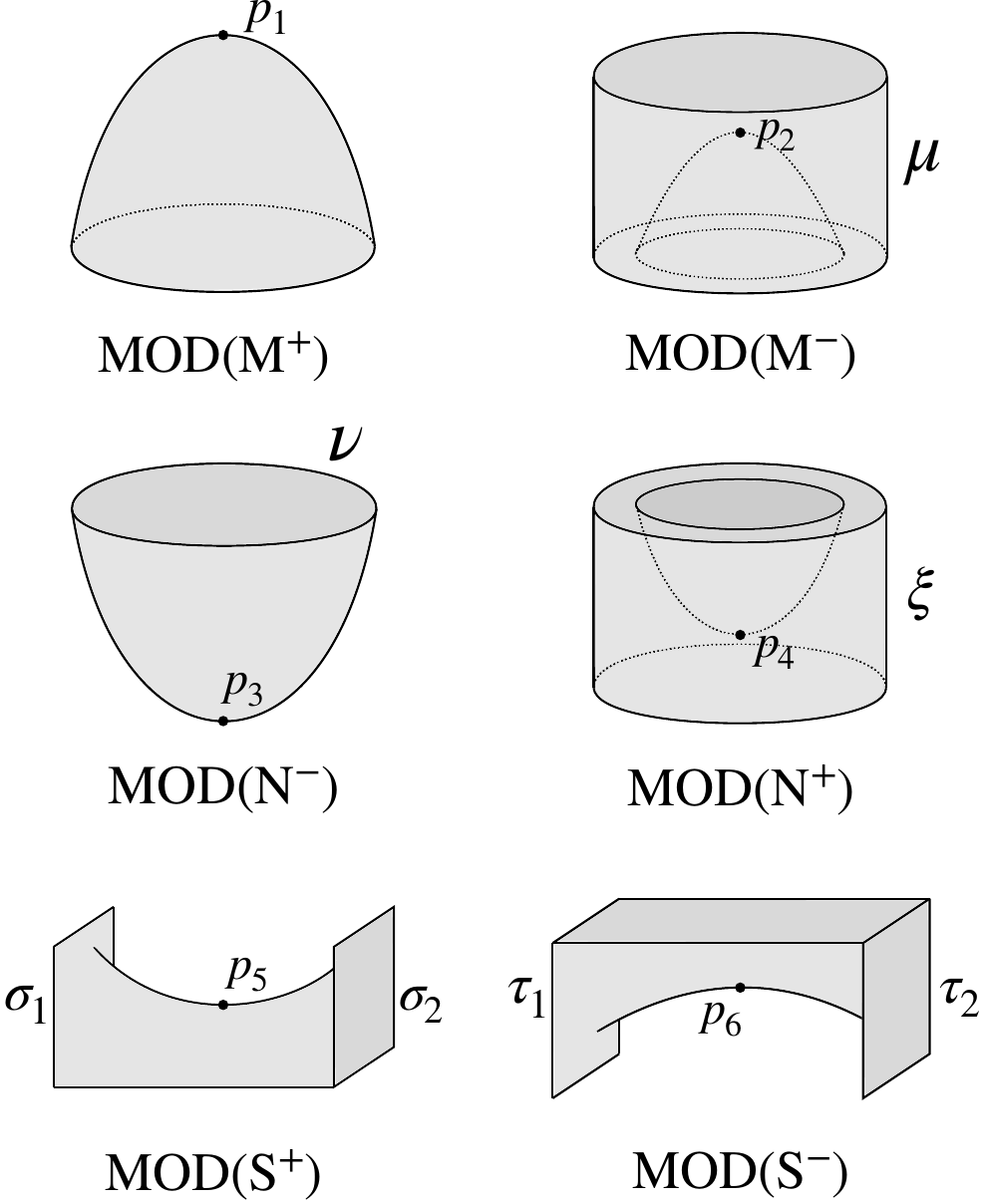}
\vspace{5pt}
\caption{
We depict the parts used to construct a $3$-dimensional manifold $N$.
The points $p_{j}$ $(j=1,2,\ldots,6)$ lie on these parts.
The cylindrical sides of $\MOD(M^{-})$ and $\MOD(N^{+})$ are represented by $\mu$ and $\xi$, respectively.
The disk side of $\MOD(N^{-})$ is denoted by $\nu$.
The rectangles on $\MOD(S^{+})$ and $\MOD(S^{-})$ are respectively denoted by $\sigma_{j}$ and $\tau_{j}$ $(j=1,2)$.
The symbols $\MOD(~)$ indicate these parts as in \cite{C}.
}
\label{fig:8}
\end{figure}

Below, we will explain how we obtain $N_{i+1}$ and $G_{i+1}$ from $N_{i}$ and $G_{i}$ to satisfy the above conditions. 
Although the methods are similar to the Curley's method in \cite{C}, we show them for readers' convenience.
Let us denote $\St(x)$ as the star of a vertex $x$ in the graph component containing $x$.

Consider the case where $C$ maps $v_{i+1}$ to $w_{i+1}$ according to $M^{+}$ in \Cref{fig:2}.
Let $\varphi_{1}\colon L_{i+1}\times\{0,1\}\to L_{i+1}\subset N_{i}$ be a projection onto the first component and $I=[0,1]$. 
We define $N_{i+1}$ as
\begin{equation*}
N_{i}\cup_{\varphi_{1}} L_{i+1}\times I\sqcup \MOD(M^{+}).
\end{equation*} 
Then, a submersion $G_{i+1}\colon N_{i+1}\to S^{1}$ is derived from $h$, where the restriction of $G_{i+1}$ to its boundary has $p_{1}\in\partial N_{i+1}$ as a critical point.

Consider the case where $C$ maps $v_{i+1}$ to $w_{i+1}$ according to $M^{-}$ in \Cref{fig:2}. 
Let $\varphi_{2}\colon\mu\to (K\setminus\Int D^{2})\times I $ be the map that sends the cylindrical side $\mu$ of $\MOD(M^{-})$ to the side formed by removing $\Int D^{2}$ from $K$, where $K$ is the connected component of $L_{i+1}$ such that an edge in $\St(w_{i+1})$ contains its corresponding point in $V$. 
We define $N_{i+1}$ as 
\begin{equation*}
N_{i}\cup_{\varphi_{1}}\Big(\overline{(L_{i+1}\setminus D^{2})}\times I\cup_{\varphi_{2}}\MOD(M^{-})\Big),
\end{equation*}
where $\varphi_{1}$ is defined similarly to the case of $M^{+}$.
Then, a submersion $G_{i+1}\colon N_{i+1}\to S^{1}$ is obtained from $h$, where the restriction of $G_{i+1}$ to its boundary has $p_{2}\in\partial N_{i+1}$ as a critical point.

\begin{figure}[t]
\centering
\includegraphics[width=115mm]{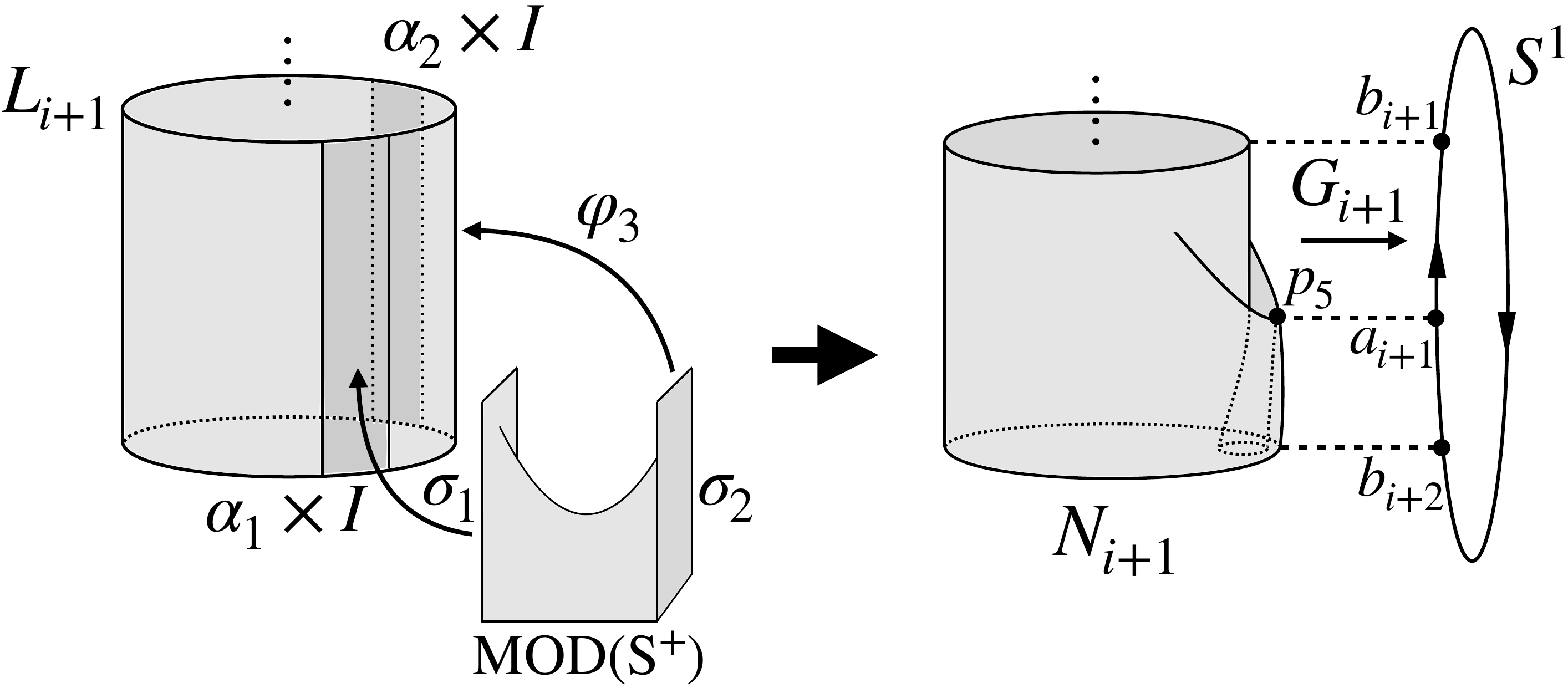}
\vspace{-3pt}
\caption{
The construction in the case of $S^{+}$.
}
\vspace{-3pt}
\end{figure}

Consider the case where $C$ maps $v_{i+1}$ to $w_{i+1}$ according to $S^{+}$ in \Cref{fig:2}. 
Let $P$ denote the connected component of $\partial L_{i+1}$ such that an edge in $\St(v_{i+1})$ contains its corresponding point in $W_{f}^{\pm}$, and let $\alpha_{1}$ and $\alpha_{2}$ be disjoint simple arcs in $P$. 
Let $\varphi_{3}\colon \sigma_{1}\sqcup \sigma_{2}\to L_{i+1}\times I$ be the map that sends the each rectangle $\sigma_{j}$ in $\MOD(S^{+})$ to $\alpha_{j}\times I$ $(j=1,2)$. 
We define $N_{i+1}$ as 
\begin{equation*}
N_{i}\cup_{\varphi_{1}}L_{i+1}\times I\cup_{\varphi_{3}} \MOD(S^{+}),
\end{equation*}
where $\varphi_{1}$ is defined similarly to the case of $M^{+}$. 
Then, we obtain a submersion $G_{i+1}\colon N_{i+1}\to S^{1}$ along $h$, where $p_{5}\in\partial N_{i+1}$ is a critical point of the map on the boundary.
See \textrm{Figure~8} for more details.

Consider the case where $C$ maps $v_{i+1}$ to $w_{i+1}$ according to $G^{+}$ or $J^{+}$ in \Cref{fig:2}. 
Let $P_{1}$ and $P_{2}$ denote the connected components of $\partial L_{i+1}$ such that two edges in $\St(v_{i+1})$ contain their corresponding points in $W_{f}^{\pm}$.
Let $\alpha_{1}$ and $\alpha_{2}$ be simple arcs in $P_{1}$ and $P_{2}$ respectively. 
We define $N_{i+1}$ as 
\begin{equation*}
N_{i}\cup_{\varphi_{1}}L_{i+1}\times I\cup_{\varphi_{3}} \MOD(S^{+}), 
\end{equation*}
where $\varphi_{1}$ and $\varphi_{3}$ are defined similarly to the case of $S^{+}$. 
Then, we derive a submersion $G_{i+1}\colon N_{i+1}\to S^{1}$ along $h$, where $p_{5}\in\partial N_{i+1}$ is a critical point of the map on the boundary.

Consider the case where $C$ maps $v_{i+1}$ to $w_{i+1}$ according to $N^{+}$ in \Cref{fig:2}. 
Let $P$ and $K$ denote the connected components of $\partial L_{i+1}$ and $L_{i+1}$ respectively, such that an edge in $\St(v_{i+1})$ contains the corresponding point in $W_{f}^{\pm}$, and an edge in $\St(w_{i+1})$ contains that in $V$.
Let $A\subset K$ be an annulus such that $P$ is a connected component of $\partial A$.
The map $\varphi_{4}\colon\xi\to \overline{(L_{i+1}\setminus A)}\times I$ sends the cylindrical side $\xi$ of $\MOD(N^{+})$ to the side formed by removing $\Int A\times I$ from $L_{i+1}\times I$. 
We define $N_{i+1}$ as 
\begin{equation*}
N_{i}\cup_{\varphi_{1}}\Big(\overline{(L_{i+1}\setminus A)}\times I \cup_{\varphi_{4}}\MOD(N^{+})\Big),
\end{equation*}
where $\varphi_{1}$ is defined similarly to the case of $M^{+}$.
This construction yields a submersion $G_{i+1}\colon N_{i+1}\to S^{1}$ along $h$, where $p_{4}\in\partial N_{i+1}$ is a critical point of the map on the boundary.

\begin{figure}[t]
\centering
\includegraphics[width=90mm]{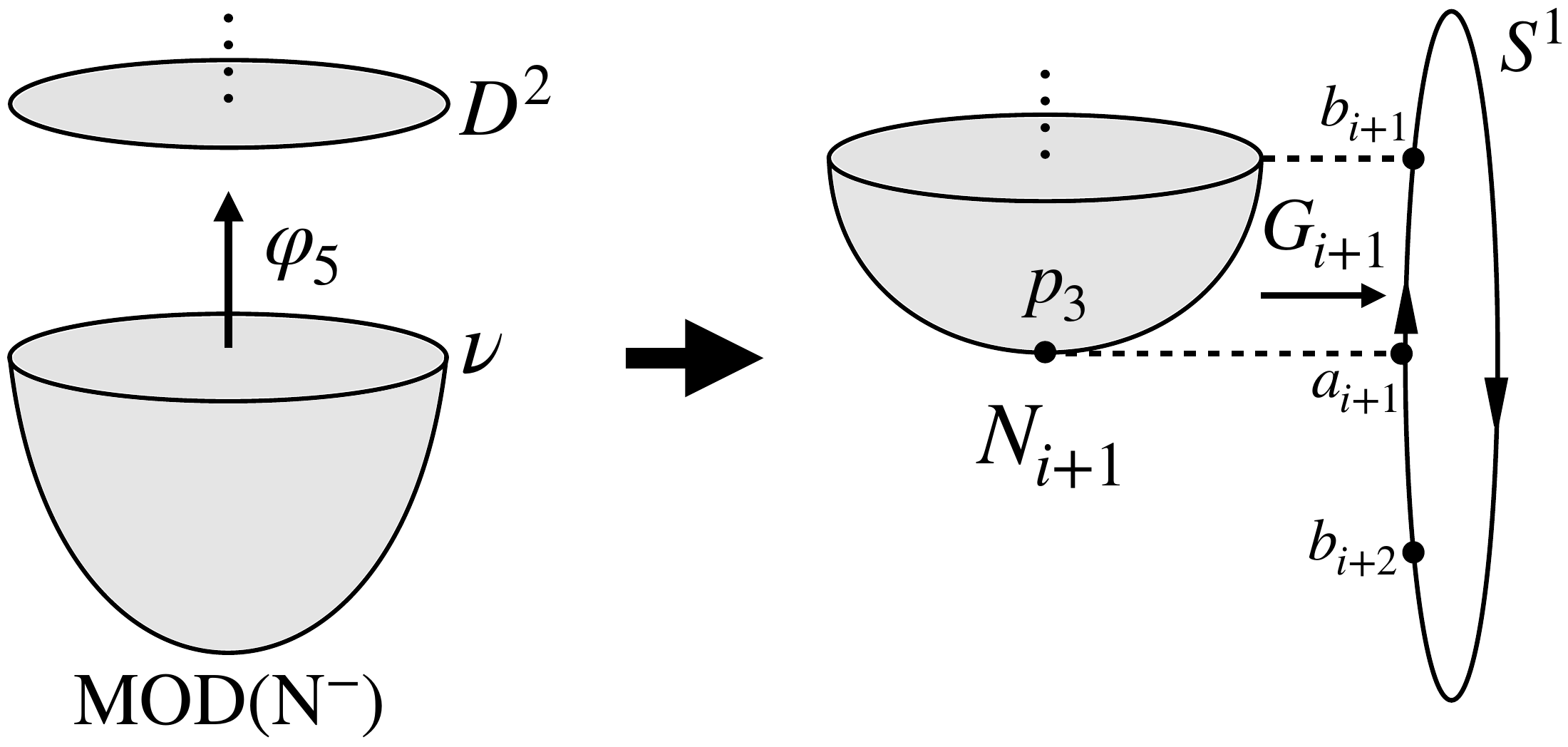}
\vspace{0pt}
\caption{
The construction in the case of $N^{-}$.
}
\vspace{0pt}
\end{figure}

Consider the case where $C$ maps $v_{i+1}$ to $w_{i+1}$ according to $N^{-}$ in \Cref{fig:2}. 
The disk $D^{2}$ is a connected component of $L_{i+1}$ such that an edge in $\St(w_{i+1})$ contains its corresponding point in $V$. 
We define $N_{i+1}$ as 
\begin{equation*}
N_{i}\cup_{\varphi_{1}\cup\varphi_{5}}((L_{i+1}\setminus D^{2})\times I\sqcup\MOD(N^{-})),
\end{equation*}
where $\varphi_{1}$ is defined similarly to the case of $M^{+}$ and $\varphi_{5}\colon\nu\to N_{i}$ attaches the disk side $\nu$ of $\MOD(N^{-})$ to $D^{2}\subset L_{i+1}$. 
Then, a submersion $G_{i+1}\colon N_{i+1}\to S^{1}$ is obtained from $h$, where $p_{3}\in\partial N_{i+1}$ is a critical point of the map on the boundary.
See \textrm{Figure~9} for more specific details.

Consider the case where $C$ maps $v_{i+1}$ to $w_{i+1}$ according to $S^{-}$ in \Cref{fig:2}. 
Let $P_{1}$ and $P_{2}$ denote the boundary components of $L_{i+1}$ such that two edges in $\St(v_{i+1})$ contain their corresponding points in $W_{f}^{\pm}$ respectively.
Let $K$ represent the connected component of $L_{i+1}$ such that an edge in $\St(w_{i+1})$ contains that in $V$. 
Let $R$ be a rectangle in $K$ such that it joints $P_{1}$ and $P_{2}$, 
where $\beta_{1}$ and $\beta_{2}$ denote the different sides of $R$ that are not in $P_{1}$ or $P_{2}$.
We define $N_{i+1}$ as 
\begin{equation*}
N_{i}\cup_{\varphi_{1}}\Big(\overline{(L_{i+1}\setminus R)}\times I\cup_{\varphi_{3}}\MOD(S^{-})\Big),
\end{equation*}
where $\varphi_{3}\colon\tau_{1}\sqcup \tau_{2}\to\overline{(L_{i+1}\setminus R)}\times I$ attaches each rectangle $\tau_{j}$ of $\MOD(S^{-})$ to $\beta_{j}\times I$ $(j=1,2)$, and $\varphi_{1}$ is defined similarly to the case of $M^{+}$.
Then, we derive a submersion $G_{i+1}\colon N_{i+1}\to S^{1}$ along $h$, where $p_{6}\in\partial N_{i+1}$ is a critical point of the map on the boundary.

Consider the case where $C$ maps $v_{i+1}$ to $w_{i+1}$ according to $G^{-}$ in \Cref{fig:2}.
Let $P$ and $K$ be respectively the connected components of $\partial L_{i+1}$ and $L_{i+1}$, such that an edge in $\St(v_{i+1})$ contains its corresponding point in $W_{f}^{\pm}$, and an edge in $\St(w_{i+1})$ contains that in $V$.
Let $R\subset K$ be a rectangle such that it joins $P$ itself and $K\setminus R$ is connected. 
We define $N_{i+1}$ as 
\begin{equation*}
N_{i}\cup_{\varphi_{1}}\Big(\overline{(L_{i+1}\setminus R)}\times I\cup_{\varphi_{3}}\MOD(S^{-})\Big),
\end{equation*}
where $\varphi_{1}$ and $\varphi_{3}$ are defined similarly to the case of $S^{-}$. 
This construction yields a submersion $G_{i+1}\colon N_{i+1}\to S^{1}$ along $h$, where $p_{6}\in\partial N_{i+1}$ is a critical point of the map on the boundary.

Consider the case where $C$ maps $v_{i+1}$ to $w_{i+1}$ according to $J^{-}$ in \Cref{fig:2}. 
Let $P$ and $K$ be respectively connected components of $\partial L_{i+1}$ and $L_{i+1}$ with similar properties as in the case of $G^{-}$. 
Let $R\subset K$ be a rectangle such that it joins $P$ itself and divides $K$ into two regions $K_{1}$ and $K_{2}$.
Note that the genera of $K_{1}$ and $K_{2}$ are equal to the numbers of the paths on the corresponding edges in $V$.
We define $N_{i+1}$ as 
\begin{equation*}
N_{i}\cup_{\varphi_{1}}\Big(\overline{(L_{i+1}\setminus R)}\times I\cup_{\varphi_{3}}\MOD(S^{-})\Big),
\end{equation*}
where $\varphi_{1}$ and $\varphi_{3}$ are defined similarly to the case of $S^{-}$. 
Then, we obtain a submersion $G_{i+1}\colon N_{i+1}\to S^{1}$ along $h$, where $p_{6}\in\partial N_{i+1}$ is a critical point of the map on the boundary.

We obtain $N_{n}$ and $G_{n}$ through iteration of these methods.
Subsequently, $N$ and $G$ are derived by attaching the initial surfaces $L_{1}$ to the final surfaces $L_{n+1}$.
This is done so that the boundary components of each connected component of the level surfaces of $G$ correspond to the elements of the inverse image by $C$ of the points in $V$ that correspond by $h$.
\qed

\begin{rem}
In the definition of a collapse in \cite{C}, surjectivity is required.
However, neither \cite{C}, nor the proof of \textrm{Theorem~1.1} of this paper is this property necessary. 
Moreover, a collapse induced by a non-singular extension of a submersion is not generally surjective. 
The non-singular extension $G\colon (S^{2}\times S^{1})\setminus D^{3}\to S^{1}$ of a submersion $g\colon S^{2}\times[0,1)\to S^{1}$ in \textrm{Figure~10} induces a collapse that is not surjective. 
\begin{figure}[t]
\centering
\includegraphics[width=85mm]{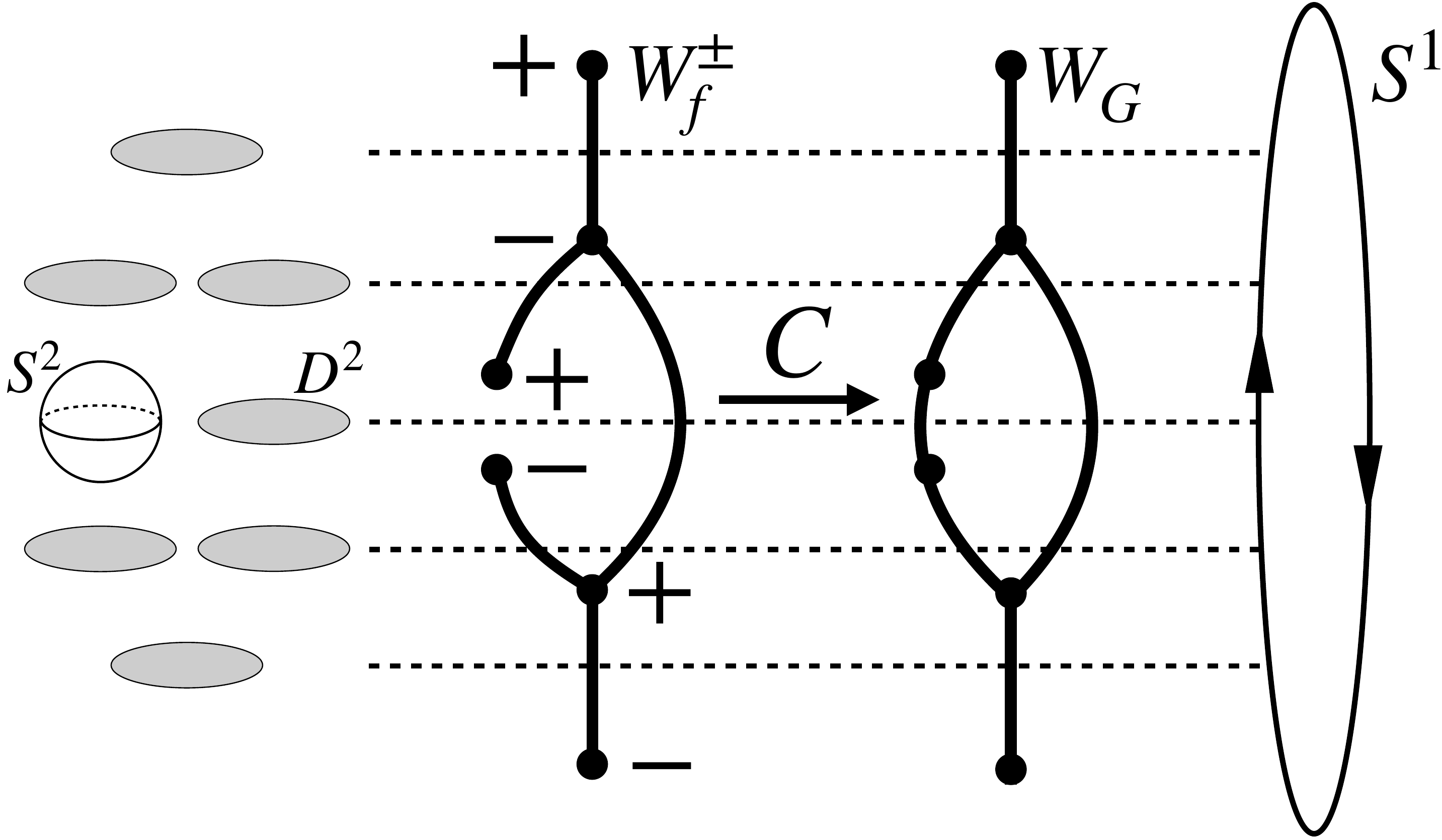}
\caption{
We depict the collapse induced by the submersion $G$ in \textrm{Remark~4.1}.
The dotted lines represent the images of $\overline{f}\colon W_{f}^{\pm}\to S^{1}$ and $\overline{G}\colon W_{G}\to S^{1}$, where $f=g|_{S^{2}\times\{0\}}$.
The surfaces on the left correspond to the level surfaces of $G$ associated with the dotted lines.
}
\vspace{-6pt}
\end{figure}
\end{rem}

\vspace{-10pt}
\section*{Aknowledgements}
The author would like to thank Osamu Saeki for fruitful comments and suggestions.
The author would also like to thank Naoki Kitazawa and Ryosuke Ota for useful advice and comments.
We appreciate the referee for careful reading and helpful comments. 
This work has been partially supported by JSPS KAKENHI Grant Number JP23H05437 and by WISE program (MEXT) at Kyushu University.

\end{document}